\documentclass [10pt] {article}
\usepackage{amsthm, amsmath}
\usepackage[all]{xy}
\usepackage{latexsym}
\usepackage{verbatim}
\usepackage{enumerate}
\usepackage{amsthm}
\usepackage{amsmath}
\title{Flops and  Equivalences of derived Categories for   Threefolds with only terminal  Gorenstein Singularities}
\author{Jiun-Cheng Chen}

\begin {document}
\maketitle
\baselineskip 12pt

\newcounter{theorem}[section]
\numberwithin{equation}{section}

\newtheorem{cl}[theorem]{Claim}

\newtheorem{Thm}[theorem]{Theorem}
{\theoremstyle{remark}
 
\newtheorem{Rem}[theorem]{\text{\textbf{Remark}}} }
\newtheorem{Def}[theorem]{Definition}
\newtheorem{Lem}[theorem]{Lemma}
\newtheorem{Prop}[theorem]{Proposition}
\newtheorem{Cor}[theorem]{Corollary}
\newtheorem{Ex}[theorem]{Example} 
{\theoremstyle{definition}
 \newtheorem{Remark}[theorem]{Remark}} 

\newcommand{\bL}{\mathop{\otimes}\limits^{L}}
\newcommand{\mExt}{\mathcal{E}xt}
\newcommand{\mHom}{\mathcal{H}om}
\newcommand{\R}{\mathbf{R}}
\newcommand{\cL}{\mathbf{L}}
\newcommand{\mO}{\mathcal{O}}  
\newcommand{\mX}{\mathcal{X}}
\newcommand{\mY}{\mathcal{Y}}
\newcommand{\mW}{\mathcal{W}}
\newcommand{\mS}{\mathcal{S}}
\newcommand{\mT}{\mathcal{T}}
\newcommand{\mC}{\mathcal{C}}
\newcommand{\Hom}{\mathrm{Hom}}

\newcommand{\mA}{\mathcal{A}}

\newcommand{\Ext}{\mathrm{Ext}}
\begin{abstract}
The main purpose of this paper is to show that  Bridgeland's moduli space of perverse point sheaves  for certain flopping contractions gives the flops, and the Fourier-Mukai transform given by the birational correspondence of the flop is an equivalence between bounded derived categories.  
\end{abstract}
 
\section{Introduction}
\subsection{The minimal model program}
One of the most important problems in birational geometry is the minimal model program (MMP). The main goal of the MMP is   to find in each birational class of varieties some distinguished representatives (minimal models)   which are ``easier'' to understand, then to use these minimal models to study the birational properties of varieties. 
In dimension 2,  satisfactory  answers have been  known for a long time.   The procedure for producing a minimal model for $X$ is repeatedly contracting a $(-1)$-curve.  The final result  of the MMP for a non-ruled surface is a smooth surface  such that it  is minimal in the category of smooth surfaces (minimal in the classical sense), and its canonical bundle  is {\it nef} (minimal in the sense of the MMP).  In higher dimensions, the situation is much more complicated. 
Certain kinds of singularities are needed  even if we start with a smooth variety. Besides singularities, we also need to consider  flops and flips, which do not occur in dimension 2. 
Based on contributions from Reid, Mori, Kawamata, Koll\'{a}r, Shokurov and others, the MMP  program was completed in dimension 3  by Mori in 1988. 

The proof of  the MMP in dimension 3  uses a very careful analysis on 
two-dimensional Du Val singularities and threefold singularities. It is very difficult to generalize the proof along these lines to higher dimensions. A more conceptual proof is very desirable.

   Flops  can be considered as a  sort of ``birational surgery'', an analogue of   surgery in algebraic topology. 
A  very natural and interesting question is what kind of invariants remain the same under flops. 
An example in this direction is that two birational nonsingular  Calabi-Yau manifolds  have the same Hodge numbers (see [Ba97] for a result on Betti numbers, or [Wa98] for a more general theorem). In dimension 3,  this theorem was first proven in [Ko89].
\subsection{Flops and derived categories}
Following Bondal-Orlov \cite{bo:gnus} and Bridgeland \cite{br:gnus}, it is plausible that the  MMP may be understood in the context of derived categories. Given a variety $X$, the minimal model(s) might be viewed as some minimal triangulated subcategories inside $D^{b}(X)$. In this picture,  it is very natural to view flops as taking different triangulated  sub-categories which are equivalent to one another, and flips as taking suitable fully faithful triangulated subcategories. There is considerable  evidence to support this picture. A very important and  interesting theorem to support this picture  is a theorem by Bridgeland.
  
In \cite{br:gnus} Bridgeland gives a moduli construction of smooth threefold flops. The moduli space he constructs is actually a 
fine moduli space. 
Furthermore, he is able to prove   a result on the equivalence of derived categories    by using techniques in [Br98] and \cite{bkr:gnus}. As a corollary of his theorem, he proves again that two birational nonsingular  Calabi-Yau threefolds have the same Hodge numbers. An interesting question is: Is  Bridgeland's theorem true for singular varieties? In this paper we  generalize his theorem to 
threefolds  with  terminal Gorenstein singularities. 
We remark here that these  singularities are isolated hypersurface singularities (see \cite{km:gnus} p.169). 
The main theorem in our paper is: 
\begin{Thm}\label{main1}
Let $X$ be a quasi-projective threefold with only terminal Gorenstein singularities and let  $f: X \rightarrow Y$ be a flopping contraction. Denote by $X^{+}$ the flop. Denote by $W=W(X/Y)$ the distinguished component of the moduli space of perverse point sheaves $M(X/Y)$ {\rm(see Appendix A for the definitions).} Then

\begin{enumerate}[{(}1{)}]
\item $W$ has only terminal Gorenstein singularities, 

\item the Fourier-Mukai type  transform $\Psi: D^{b}(W) \rightarrow D^{b}(X)$ induced by the universal perverse point sheaves is an equivalence, and

\item $W \cong X^{+}$. 

\end{enumerate}   

\end{Thm}

\subsection{Reduction to the local cases}
We   outline the proof of this theorem in the subsequent subsections. First, a few  comments on Fourier-Mukai type transforms. A  Fourier-Mukai type transform F  may not send $D^{b}(W)$ to $D^{b}(X)$ 
since $X$ and $W$ may be singular.  However, the kernel we consider is $[ \; \mathcal{I} \rightarrow \mO_{W \times X}\;]$, where $\mathcal{I}$ is  the universal perverse ideal sheaf and hence is flat over $W$. We show in Section 2 
that such a kernel does define a transform $\Psi: D^{b}(W) \rightarrow D^{b}(X)$. Let $\{Y_{i}\}$ be an affine cover of $Y$.
We pull back this universal perverse point sheaf to each $Y_{i}$.
These kernels  give Fourier-Mukai type transforms 
$\Psi_{i}: D^{b}(W_{i}) \rightarrow D^{b}(X_{i})$. 

We note that Theorem~\ref{main1} is local in $Y$. 
Since the moduli space $W$ is local in $Y$ (see Remark~\ref{quai} in Appendix A),  part (1) and  part (3) of  Theorem~\ref{main1} are  clear. 
It is not obvious that part (2) is also local in $Y$ since we can not check whether a functor is an equivalence  or not locally.
The next proposition shows that  part (2)
 of the theorem is also local
in $Y$.  
The main point of the proof  is that $\Psi$ has a right adjoint $\Phi$. 
\begin{Prop}{\rm (see Proposition~\ref{fin})}\label{fin0}
Notation as above. If there is an affine cover $\{ Y_{i} \}$ of $Y$ such that  $\Psi_{i}: D^{b}(W_{i}) \rightarrow D^{b}(X_{i})$ are equivalences of derived categories, then $\Psi:D^{b}(W) \rightarrow D^{b}(X)$ is also an equivalence of derived categories.
\end{Prop}   
\subsection{Results from \cite{br:gnus} and \cite{bkr:gnus}}
The argument in \cite{br:gnus} uses the 
 non-singularity   assumption  in a significant way.
 The techniques used in his proof do not 
  seem to generalize directly to singular varieties. Our idea  is that instead of studying the singular threefold directly, we study a nonsingular fourfold, which is a smoothing,  and see how much information  about the singular threefold we can  get from this smooth fourfold. 
The starting point of  our approach in this paper is   
the following theorem, which   is  a restatement of a  combination of  results in [Br00] and [BKR99]. We sketch  a proof in Appendix~\ref{sec:2ndapp} for the reader's convenience.  
\begin{Thm}\label{bkr}{\rm (\cite{br:gnus} and  \cite{bkr:gnus})} 
Let $X \rightarrow Y$ be a flopping contraction where $X$ is an $n$-dimensional smooth quasi-projective variety and the dimension of every fiber is at most $1$. Let $W$ be the distinguished component of the moduli space of perverse point sheaves $M(X/Y)$. Assume \nolinebreak{${\rm dim} (W \times _{Y} W) \leq n+1$}. 
Then
\begin{enumerate}[{(}1{)}]
\item $W$ is smooth, 

\item the transform $D^{b}(W) \rightarrow D^{b}(X)$ induced by the universal perverse point sheaf is an equivalence, and 

\item The flop $X^{+}$ for $ X \rightarrow Y$ exists and $W \cong X^{+}$.

\end{enumerate}

\end{Thm}

The next corollary follows immediately from Theorem~\ref{bkr}.  

\begin{Cor}\label{bkr1}
Let $f:X \rightarrow Y$ be a flopping contraction with $X$ smooth.  Then
\begin{enumerate}[{(}1{)}]
\item if ${\rm dim} (X)=3$, then  the conclusions in Theorem~\ref{bkr}  always hold, and
  
\item if ${\rm dim} (X)=4$, every fiber of  $f: X \rightarrow Y$ is of dimension at most $1$, and $g: W \rightarrow Y$  does not contract any divisor to a point, then  the conclusions in Theorem~\ref{bkr}  hold. 
\end{enumerate}
\end{Cor}

The 
next proposition  is a combination of results in [BKR99] and [Br00] as indicated by Bridgeland in the introduction in [Br00]. We shall not need this result in this paper.
\begin{Prop}{\rm (see [Br00] and  [BKR99])} \label{irr}
Notation as in Theorem~\ref{bkr}. Assume that ${\rm dim}(X)=3$. Then $M(X/Y)=W(X/Y)$. 
\end{Prop}  
\subsection{Relations between $W(X/Y)$ and $W(S/T)$ for a Cartier divisor $T \subset Y$} 
Let $X$ be a variety with at worst terminal Gorenstein  singularities. Assume that ${\rm dim}(X)$ is either  $3$ or  $4$ and $f:X \rightarrow Y$  satisfies the  following two conditions
\begin{enumerate}[{(B.}1{)}]\label{con}
\item $\R  f_{*} \mathcal{O}_{X} = \mathcal{O}_{Y}$, and
\item   every  fiber of $f$ is of ${\rm dim} \leq 1$.
\end{enumerate}
Let $T \subset Y$ be an effective Cartier  divisor; for simplicity we assume that it is  an integral subscheme of $Y$.
Let $S$ be the preimage of $T$ in $X$.  Denote by $W_{T} =W (X/Y)_{T}$ the restriction of the moduli space $W(X/Y)$ to  $T \subset Y$.  
The underlying philosophy of our approach is that
\begin{enumerate}[{(}1{)}]
 \item  we find a smoothing $F: \mX \rightarrow \mY$ of  $f: X \rightarrow Y$, and 
  
 \item  we  relate the fiber  of the  moduli spaces $W(\mX/\mY)$ to the moduli space $W(X/Y)$.
\end{enumerate}

The next proposition shows that $(2)$ is possible. 
\begin{Prop}{\rm (see Proposition~\ref{dim3})}\label{goingdown1}
Notation  as above. There is a canonical morphism $W(X/Y)_{T} \hookrightarrow M(S/T)$, which is an inclusion of components.

\end{Prop}

\begin{Remark}
It is also true that $M(X/Y)_{T}=M(S/T).$
We shall not need this stronger result in our paper. 
\end{Remark}

\subsection{Smoothing and smooth hyperplane sections}
The  following proposition shows that smoothing is always possible after passing to an affine cover.
\begin{Prop}{\rm (see Proposition~\ref{smo} for the precise statement)} \label{smo1}
Let $X \rightarrow Y$ be a flopping contraction between threefolds where $X$ has at worst terminal Gorenstein singularities. Then there is an affine cover $\{ Y_{i} \}$ of $Y$ such that each  $f_{i}: X_{i} \rightarrow Y_{i}$ is   a  smoothable morphism.
\end{Prop}
In the remainder of this subsection and the next subsection, we work over $Y_{i}$. We shall suppress the indices  when no confusion is possible. 

Let $F: \mX \rightarrow \mY$ be a  one-parameter  deformation of $f: X \rightarrow Y$ such that $\mX$ is nonsingular. Let $Y_{sing}= \{ p_{i} : i=1, \cdots , m \}$ be the finite set of singular points of $Y.$ We also consider  them as points of $\mY.$  Let $\mathcal{T}$ be a general hyperplane section passing through  $Y_{sing} \subset Y \subset \mY.$  Denote by  $\mathcal{S}$ the preimage of $\mathcal{T}$. 

The following  proposition enables us to use results on  the smooth  threefolds in [Br00]. 
\begin{Prop}{\rm (see Proposition~\ref{gen} for the precise statement)} \label{gen1}
The hyperplane section $\mathcal{S}$ is nonsingular. 
\end{Prop}
\subsection{Proofs of  ``no divisor is contracted to a point'' and Theorem~\ref{main1}}
Denote by $\mW$ the distinguished component of the  moduli space of perverse point sheaves for $F: \mX \rightarrow \mY.$ Denote by $G : \mW \rightarrow \mY$  the natural birational morphism.   
We   explain briefly how to prove that $G: \mW \rightarrow \mY$ contracts no divisor to a point. Assume that there is a divisor contracted to a point, say $p$, by $G$.  This point $p$ must be one of the singular points in $Y.$
Take a general hyperplane section $\mT$ of $\mY$ passing through $p $. The preimage of $\mT$, denoted  by $\mS \subset \mX$, is smooth. By \cite{br:gnus}  the connected component of $W(\mS/ \mT) \subset M(\mS/ \mT)$ is smooth. Every  component $W_{j}$ of $W(\mX/\mY)_{\mT}$ is a component of $M(\mS/\mT)$ by  Proposition~\ref{goingdown1}. The fiber $W(\mX/\mY)_{\mT}$  is connected. Therefore, the distinguished component  $W(\mS/ \mT)$ is the only component by the smoothness result. Since $W(\mS/\mT) \rightarrow \mT$ is birational,
 it follows that the preimage of $p$ is at most two-dimensional, a contradiction.   

By Corollary~\ref{bkr1}, it follows that   the Fourier-Mukai type transform $D^{b}(\mW) \rightarrow D^{b}(\mX)$ is an equivalence
and $\mW \cong \mX^{+}$.

By standard results on flops, it follows easily that $W \cong X^{+}$ is the flop and hence has only terminal Gorenstien singularities (see Section~\ref{sec:going}). This concludes the proof of part (1) and (3) in Theorem~\ref{main1}. \newline

To prove $D^{b}(W) \cong D^{b}(X)$,   more work is needed. 
Let  $\Psi: D^{b}(\mW) \rightarrow D^{b}(\mX)$ be the Fourier-Mukai type transform defined by the universal perverse point sheaf, i.e. the  structure sheaf of the fiber product $\mW \times _{\mY} \mX$. 
Let $i_{0}$ be  the inclusion morphism $W \rightarrow \mW$ (see Proposition~\ref{dim3}).
Denote
by $\Psi_{0}: D^{b}(W) \rightarrow D^{b}(X)$ the Fourier-Mukai type transform defined by $\cL i^{*}_{0} \mO_{\mW \times _{\mY} \mX}.$ Note that this 
 Fourier-Mukai type transform  is equivalent to the Fourier-Mukai type transform defined by the kernel 
$\mO_{W \times_{Y} X}$ (see Proposition~\ref{dim3} and  Corollary~\ref{universalobject}). We shall denote both of these two functors by $\Psi_{0}$.
Denote by $\Phi: D^{b}(\mX) \rightarrow D^{b}( \mW)$ the right adjoint to $\Psi$.  This functor is also a Fourier-Mukai type transform (see Lemma 4.5 in [Br98]). \newline

To complete the proof of part (2) in Theorem~\ref{main1}, we use
the next proposition. The proof of this proposition
is given in Section~\ref{sec:going}. The main point  is 
to   show  that $\Psi(i_{0\;*}(-))\cong i_{0\;*}  (\Psi_{0}(-))$.
\begin{Prop}\label{fourto3}{\rm (see Proposition~\ref{4to3})}
\[ \Psi: D^{b}(\mW) \cong D^{b}(\mX)\;  \Longrightarrow \;  \Psi_{0}: D^{b}(W) \cong D^{b}(X). \]
\end{Prop}

\begin{Remark}
Using the results in [Ne99] and [Ne00], our results  imply the $K$-theories of coherent sheaves (i.e. $G$-theories)  of $X^{+}$ and $X$ are isomorphic.
\end{Remark}
\subsection{Comments and further developments}
Finally, we would like to say a few words on the limitation of the smoothing approach and our speculation on the possible generalizations. 

It is well-known that quotient singularities in dimension $\geq 3$ are rigid. Therefore our smoothing  approach would not work for the most general threefold flops.  To settle   general  three-dimensional flops  using Bridgeland's approach, it seems that  new ideas and techniques are needed.
We speculate that algebraic  stacks should play certain roles in the complete picture. 
Recently Kawamata proved an interesting result on $n$-dimensional toric flips and derived categories (see \cite{ka:gnus}). His result  provides some evidence to support our speculation. 

In the flips cases, 
D. Abramovich and I are working on some simple toric flips ([AC01]).  In that case, we use the natural stack structure on threefolds in question  instead of using deformations.  We also  plan to use the  similar stack structure 
to extend our results to Q-Gorenstein case.

\subsection{Plan of the paper}
The plan of this paper is  as follows. 
In Section 2  we present  a few basic facts about the Fourier-Mukai 
type transforms. 
In Section 3, we explain how to reduce the proof to an affine Zariski neighborhood of $Y$. 
In Section 4, 
we prove several facts on the
moduli space of perverse point sheaves. 
In Section 5, we give the proofs of lemmas on the deformation and general hyperplane sections  needed for our proof.   We give a proof on how to deduce the equivalence of derived categories  in dimension $3$ from the corresponding result in dimension $4$ in Section 6.

The first appendix contains basic facts about triangulated categories and perverse coherent sheaves. All the material is taken from \cite{br:gnus}. 
We sketch  the proof of Theorem~\ref{bkr} in the second appendix. The proof is  the same as the proof in [BKR99].  
   
\subsection{Notation}
All schemes $T$ are schemes of finite type over \textbf{C}. Denote by $T^{n}$ the normalization of $T$. 
Denote by  $D_{qc}(T)$ the derived category of the abelian category $Qcoh(T)$ of quasi-coherent $\mathcal{O}_{T}-$modules.
Denote by $D^{+}(T)$ the full subcategory of  $D_{qc}(T)$  consisting of complexes whose cohomology sheaves  are   bounded below and  coherent. Denote by $D^{-}(T)$ the full subcategory of  $D_{qc}(T)$ consisting of complexes whose    cohomology sheaves  are   bounded  above and  coherent. 
Denote by  $D^{b}(T)$ the full  subcategory of $D_{qc}(T)$ consisting of complexes  with bounded and coherent cohomology sheaves. Denote by  $D^{b}_{c}(T)$ the full  subcategory of $D^{b}(T)$ consisting of complexes whose cohomology sheaves are of  proper support.  

Let $f: T \rightarrow S$ be a projective birational morphism such that the conditions (B.1) and (B.2)  are satisfied. We denote by $M(T/S)$ the fine moduli space of perverse point sheaves. 
Let $U \subset S$ be the maximal open set such that $f^{-1} \mid_{X_{U}}$ is an isomorphism.
 Denote by  $W$ the irreducible component of $M(T/S)$ which 
contains   $U \subset S$. 

\subsection*{Acknowledgments}

This paper would not exist without  very helpful discussions with many people. The author would like to thank Tom Bridgeland for his help. His comments on the first version of this paper are very essential to this  current version.  Without his help, the author would not be able to  prove the  result on the equivalence of derived categories.  The author would also like to express his thanks to  J. de Jong,  J. Koll\'ar, A. Neeman,  A. Polishchuk,   M. Reid, J. Starr and S.T. Yau. Last but not least,  the author would like to thank Dan Abramovich, his thesis advisor, for his guidance and encouragement.  Without his insistence, the author would not  have tried  to prove the equivalence of categories result.

\section{Fourier-Mukai type transforms on singular varieties}

This section contains several basic lemmas on Fourier-Mukai type transforms.  We essentially follow [BO95]. 
\subsection{Boundedness of a transform}
Let $X$ and $Y$ be quasi-projective varieties. Consider the diagram \newline
\xymatrix{&&& &X \times Y \ar[dl]_{p_{1}} \ar[dr]^{p_{2}} & \\ &&& X&& Y.}\newline 
\newline One  can use   the  formula  
\[  \R p_{2\;*} (\mathcal{E} \bL       \cL p_{1}^{*}(-)) \]
to define  a functor $F: D_{qc}(X) \rightarrow D_{qc}(Y)$  by results in \cite{sp88}.
When  $X$ and $Y$ are smooth, every object $\mathcal{E} \in D^{b}(X \times Y)$ is of finite Tor-dimension. If $p_{2}|_{ {\rm Supp}(\mathcal{E})}: {\rm Supp}(\mathcal{E}) \rightarrow Y$ is proper, then the functor $F$
  is  also a functor of 
triangulated categories $ F:D^{b}(X) \rightarrow D^{b}(Y).$

However, an object $\mathcal{E} \in D^{b}(X \times Y)$ may not be of finite Tor-dimension when $X$ and $Y$ are not smooth. Hence  this transform $F$  may not send $D^{b}(X)$ to $D^{b}(Y)$.    
The next easy lemma shows that  many such transforms $\R p_{2 \;*}(\mathcal{E}\bL                   \cL        p_{1}^{*}(-))$  send $D^{b}(X)$ to $D^{b}(Y)$. 
\begin{Lem}\label{finite}
Assume that $\mathcal{E} \in D^{b}(X \times Y)$ is
isomorphic to a complex F of coherent 
$\mathcal{O}_{X \times Y}$-sheaves such that
each of these sheaves is flat over $\mathcal{O}_{X}$, 
and ${\rm Supp}(\mathcal{E}) \rightarrow Y$ is a proper morphism. Then $\R p_{2 \;*}(\mathcal{E} \bL   \cL        p_{1}^{*}(-)) $ sends $D^{b}(X)$ to $D^{b}( Y)$.
\end{Lem}
\textbf{Proof}. We first
check the  functor  $\mathcal{E}\bL \cL p^{ *}_{1}(-)$ sends $D^{b}(X)$ to $D^{b}(X \times Y)$. This can be checked locally and follows from the identity:
\[ (M \bL                  _{C}(N \bL                  _{A} C)) \cong (M \bL                 _{A} N), \]
where $C$ is a ring flat over $A$ and $M$ is a finite complex of finitely      presented  $C$-modules and $N$ is a finite  complex of $A$-modules. Our assumption  on Tor-dimension (of $M$ over $A$)  implies that  $(M \bL_{A} N)$  
  is a finite complex of finitely presented
  $C$-modules  when $N$ is a finite complex of  finitely presented   $A$-modules.  

Let $\mathcal{F}$ be any object in $D^{b}(X)$.
Write 
$\mathcal{G}= \mathcal{E} \bL \cL p_{1}^{*}(\mathcal{F}) \in D^{b}(X \times Y)$. 
Note that   ${\rm Supp} (\mathcal{G}) \subset {\rm Supp} (\mathcal{E})$ is
  a closed subset, 
so  $\R p_{2 \;*}(\mathcal{G}) \in D^{b}(Y)$  
 by the assumption that ${\rm Supp}( \mathcal{E}) \rightarrow Y$ is proper. $\Box$

\begin{Lem}\label{fam}
Let $Z$ be a closed subscheme of $X \times Y$ and $\mathcal{E}$ be an object in $D^{b}(Z)$. Denote by $i : Z \rightarrow X \times Y$ the inclusion map. Then we have 
\[\R p_{2 \; *}( i_{*}\mathcal{E} \bL                   \cL        p_{1}^{*}(-)) \cong \R  (i \circ p_{2})_{*}( \mathcal{E} \bL                   \cL        (i \circ p_{1})^{*}(-)).\]

\end{Lem}
\textbf{Proof}. This follows easily from the projection formula. $\Box$ \newline

We use this lemma to prove the following fact. Consider the diagram \newline
\xymatrix{ &&&&& \mathcal{X} \times_{\mathcal{C}} \mathcal{Y} \ar[d]_i & \\ 
 &&&&& \mathcal{X} \times \mathcal{Y} \ar[ld]_{p_{1}} \ar[rd]^{p_{2}} & \\ 
 &&&&\mathcal{X} \ar[rd]_{f} & & \mathcal{Y} \ar[ld]^{g}  \\ &&&&& \mathcal{C}\; \;.& } \newline
Let $\mathcal{E} \in D^{b}(\mathcal{X} \times \mathcal{Y})$ be an object  which  comes from $D^{b}(\mathcal{X} \times_{\mathcal{C}} \mathcal{Y}).$  Then the Fourier-Mukai type transform $F_{\mathcal{E}}$ can be defined as
\[ \R (p_{2} \circ i)_{*} ( \cL        ( (p_{1} \circ i)^{*}(-) \bL                    \mathcal{E})) \]
by the lemma. 

\subsection{Compositions of Fourier-Mukai type transforms} 
The next proposition shows  that   the composition of two Fourier-Mukai type transforms is still a Fourier-Mukai type transform. This is a generalization of Proposition 1.4 in [BO95].

Let $X$, $Y$ and $Z$
be      quasi-projective      varieties and  $I$, $J$   objects of $D^{b}(X  \times Y)$ and  $D^{b}(Y \times Z)$(resp.). We assume that $I$ and $J$  satisfy  the  assumptions in Lemma~\ref{finite}.

Consider the  diagram of projections \newline
\xymatrix{&&&  X \times Y \times Z \ar[ld]_{p_{12}} \ar[d]_{p_{13}} \ar[dr]^{p_{23}}&& \\
&&X \times Y \ar[d]_{\pi_{12}^{1}} \ar[dr]_{\pi_{13}^{1}}  & X \times  Z \ar[dl]_{\pi_{12}^{2} } \ar[dr]_{\pi_{23}^{2}}  & 
Y \times Z \ar[dl]_{ \pi_{13}^{3}} \ar[d]^{ \pi_{23}^{3}} \\
&& X & Y &Z}\newline
 and the functors
 \[ F_{I}: D^{b}(X) \rightarrow D^{b}(Y), \]
 \[    F_{J}: D^{b}(Y) \rightarrow D^{b}(Z), \]
   defined by the formulas 
\[  F_{I}= \R  \pi_{12 \;*}^{2}(I \bL                   \cL         \pi_{12}^{1\;*}(-)),\]   
\[ F_{J}= \R  \pi_{23 \;*}^{3}(J \bL                   \cL         \pi_{23}^{2\;*}(-)).\]

\begin{Prop}\label{comp}
The composition functor of
 $F_{I}$ and $F_{J}$ is isomorphic to $F_{K}$ with 
  \[ K= \R p_{13 \; *}(\cL p_{23}^{*}J \bL                \cL  p_{12}^{*}I). \]

\end {Prop}
\textbf{Proof}. We follow the argument  in \cite{bo:gnus}:  
\begin{eqnarray}
F_{J} \circ F_{I} &=& \R  \pi^{3}_{23\;*}
                    (J \bL  \cL \pi^{2\;*}_{23}
                    (\R \pi^{2}_{12\;*}
                    (I \bL
                    \cL \pi^{1\;*}_{12}(-)))) \nonumber \\
                  &\cong&   \R  \pi^{3}_{23\;*}
                     (J \bL \R p_{23\;*}
                     (\cL p_{12}^{*}(I \bL 
                   \cL \pi^{1\;*}_{12}(-)))) \label{2.1} \\
                  &\cong& \R \pi^{3}_{23\;*}\R p_{23\;*}
                   (\cL p_{23}^{*}J \bL
                   (\cL p_{12}^{*}(I \bL 
                   \cL \pi^{1\;*}_{12}(-)))) \label{2.2} \\
                 &\cong&     \R \pi^{3}_{13\;*}
                   \R p_{13\;*}  (\cL p_{23}^{*}J
                   \bL (\cL p_{12}^{*}(I)  
                  \bL \cL p_{12}^{*}
                  \cL \pi_{12}^{1\;*}(-))) \label{2.3}\\
                   &\cong& \R \pi^{3}_{13\;*}
\R p_{13\;*}(( \cL p_{23}^{*}J 
\bL \cL p_{12}^{*}
I) \bL
\cL p_{13}^{*}\cL \pi_{13}^{1\;*}(-)) \label{2.4} \\
&\cong&
\R \pi^{3}_{13\;*}
(\R  p_{13\;*}(\cL p_{23}^{*}J
\bL \cL p_{12}^{*}I) 
\bL \cL\pi^{1\;*}_{13}(-)). \label{2.5}
\end{eqnarray}
The isomorphism \eqref{2.1} follows from the flat base change theorem, the isomorphisms \eqref{2.2} and \eqref{2.5} follow from the
projection formula. The isomorphisms \eqref{2.3} and \eqref{2.4} are obvious. $\Box$ 

\section{Reduction of the proof to affine cases}\label{sec:rop} 
Let $f : X \rightarrow Y$ be a flopping contraction 
between two quasi-projective three-dimensional normal varieties. Assume that  the variety $X$ has only terminal Gorenstein singularities.  
We explain in this section how to  reduce the proof of part (2) in Theorem~\ref{main1} to an affine cover $\{ Y_{i} \}$. Consider the  diagram \newline  
\xymatrix{&&& W \ar[dr]_{g} && X\ar[dl]^f\\
           &&&& Y.& }
\newline Fix an affine cover  $\{Y_{i}\}$ of $Y$. Pull back everything to ${Y_{i}}$ \newline
\xymatrix {&&  W_{i} \ar[dr]_{g_{i}}&& X_{i}\ar[dl]^{f_{i}}&  W \ar[dr]_{f} && X\ar[dl]^g &   \\
&&&Y_{i}&\ar[r] &&Y.&& &\\} 
\newline     Let $F: D^{b}(W) \rightarrow D^{b}(X)$ be a Fourier-Mukai type  transform defined by an object
 $\mathcal{E} \in D^{b}(W \times X)$. Assume that $\mathcal{E}$  is of finite  Tor-dimension over $W$ and  the projection morphism $W \times X \rightarrow X$  is proper when restricted to ${\rm Supp}(\mathcal{E}) \rightarrow X$.
  
Denote by $F_{i}$  the corresponding Fourier-Mukai type transforms when we pull back everything to $Y_{i}$.
Note that any Fourier-Mukai type transform
  also defines a functor on $D_{qc}$.
We show 
that if we can check the equivalence of categories locally, then by the existence of a  global adjoint functor, we are able to prove  the equivalence of derived categories. 
\begin{Remark}\label{unsatified}
In the proof on Lemma~\ref{rig}, we  need to work on $D_{qc}$ since we  
 invoke a theorem by Neeman on a very general form of Grothendieck duality (see \cite{Ne96:gnus}). 
Since  this is the only reason for passing to $D_{qc}$,   we would like to have a proof without using these huge categories.
For the time being, however, we are not  able to give such a proof. 
\end{Remark}
\begin{Prop}{\rm (= Proposition~\ref{fin0})}\label{fin}
Notation as above. Assume that all  $F_{i}: D^{b}(W_{i}) \rightarrow D^{b}(X_{i})$ are equivalences of derived categories.
Then  the Fourier-Mukai type transform  $F: D^{b}(W) \rightarrow D^{b}(X)$ is an equivalence  of derived categories.
\end{Prop}

We give several lemmas needed for the 
proof in the subsequent  subsections. The proof of this proposition is given at the end of this section.

\subsection{A spanning class}
We recall the definition of spanning classes for a triangulated category $\mA$ (see Definition 2.1. in [Br98]).  
\begin{Def}
A subclass $\Omega$ of objects
of $\mA$ is called a spanning class for $\mA$,
 if for every object $a \in \mA$ 
\[   \Hom^{i}_{ \mA}( b, a)=0 \; \; \forall \;   b \in \Omega  \; \;
  \forall i \in Z \Longrightarrow a \cong 0, \]
\[   \Hom^{i}_{\mA }(a , b)=0 \; \; \forall  \;  b \in \Omega  \; \;
  \forall i \in Z \Longrightarrow a \cong 0. \]
\end{Def}
\begin{Lem}\label{spa}
Let $X$ be a normal projective variety with only isolated singular points $\{x_{i} : i=1, \cdots , k \}.$ Let $\Omega_{1}=\{\mathcal{O}_{x}: \; x \in X\}$ and $ \Omega_{2}=\{ \mathcal{O}_{Z}: \; {\rm Supp}(Z) \subset X_{sing}= \{ x_{i}: i=1, \cdots , k \}\} $. Then $\Omega =   \Omega_{1} \bigcup \Omega_{2}$ is a spanning class for $D^{b}(X)$.
\end{Lem}
\textbf{Proof}. (a)  We check the condition
\[   \Hom^{i}_{ D^{b}(X) }( a, b)=0 \; \; \forall \;  b \in \Omega  \; \;
  \forall i \in Z \Longrightarrow a \cong 0 \]
by   using the argument in [Br98]. 
For any object $a \in D^{b}(X)$ and any $x \in X$, there is a spectral sequence 
\[ E_{2}^{p,q} = \Ext_{X}^{p}(H^{-q}(a), \mathcal{O}_{x}) \Rightarrow \Hom_{D^{b}(X)}^{p+q}(a,  \mathcal{O}_{x}). \] \
If $a$ is non-zero, let $q_{0}$ be the maximal value of $q$ such that $H^{q}$ is non-zero.  Take any  point $x \in {\rm Supp}(a)$. There is a non-zero element of $E_{2}^{0,-q_{0}}$, which survives at the $E_{\infty}$ stage. This  gives an element of $\Hom_{D^{b}(X)}(a, \mathcal{O}_{x})$, a contradiction. 

(b) The condition 
\[   \Hom^{i}_{D^{b}(X) }(b,a )=0 \; \; \forall  \;  b \in \Omega  \; \;
  \forall i \in Z \Longrightarrow a \cong 0 \]
 is equivalent to the following statement:  
\[ a \not\cong 0 \Longrightarrow \Hom_{D^{b}(X)}(b,a) \neq 0 \;  \text{for some} \; b \in \Omega. \] 
We use a similar spectral sequence
\[ E_{2}^{p,q} = \Ext_{X}^{p}(b, H^{q}(a)) \Rightarrow \Hom_{D^{b}(X)}^{p+q}(b, a ) \] \ 
to prove this statement. 

Fix any $x \in X_{reg}$.  
\begin{cl}\label{hom}
$\Hom^{i}_{D^{b}(X)}(\mathcal{O}_{x}, a)=0 \; \; \forall \; i \; \Longrightarrow 
 x \not\in {\rm Supp}(a).$ 
\end{cl} 
It is clear that for each $i$ the sheaf  $\mHom^{i}_{D^{b}(X)}(\mathcal{O}_{x}, a)$ is a coherent $ \mO_{x}$-sheaf.   Take an affine neighborhood $U=Spec(A)$ of $x$.  There is no higher derived functor for $\Gamma(Spec \; (A) \; ,-)$. Thus $\Hom_{D^{b}(X)}^{i} \cong \mHom_{A}^{i}$.   Since $x$ is a non-singular point, the sheaf $\mathcal{O}_{x}$ has a finite flat resolution.  
Thus $\R \Hom_{D^{b}(X)}(\mathcal{O}_{x}, a) \bL  \mathcal{O}_{x} \cong \R \Hom_{D^{b}(X)}(\mathcal{O}_{x}, a \bL \mathcal{O}_{x}) $. By assumption we have $\R \Hom_{D^{b}(X)}(\mathcal{O}_{x}, a)=0$, and hence  $\R \Hom_{D^{b}(X)}(\mathcal{O}_{x}, a \bL \mathcal{O}_{x}) =0.$
Replacing $a$ by $a \bL \mathcal{O}_{x}$,  we may assume that $a$ is with proper support. 

Since $X$ is projective, $X_{reg}$ is quasi-projective. Both $a$ and $\mathcal{O}_{x}$ are with proper supports. Serre duality implies that $\Hom^{i}_{D^{b}(X)}(a, \mathcal{O}_{x})=0.$ By the argument in $(a)$ above, it follows that $x \not\in {\rm Supp}(a)$. This shows that ${\rm Supp}(a) \subset X_{sing}$, which is equivalent to Claim~\ref{hom}. \newline

Since $a \in D^{b}(X)$, there is a subscheme structure $z$ on $x_{0}$ such that  $a \in D^{b}(z)$ (i.e. every cohomology group is an $\mathcal{O}_{z}$-module).  Let $q_{1}$ be the minimal value of $q$ such that $H^{q}(a) \neq 0$. It is  clear that $\R \Hom_{D^{b}(X)}( \mathcal{O}_{z}, H^{q_{1}}(a)) \neq 0$ and its elements  survive at the $E_{\infty}$ level. This concludes the proof. $\Box$

\subsection{Right adjoints}
  
\begin{Lem}\label{rig}
Let $X$ and $Y$ be projective Gorenstein varieties and $\mathcal{E}$ an object of $D_{qc}(X \times Y).$  
Consider the diagram  \newline
\xymatrix{&&&&X \times Y \ar[ld]_{p_{1}} \ar[rd]^{p_{2}} & \\
         &&& X && Y.}\ 
\newline Denote by $F: D_{qc}(X) \rightarrow D_{qc}(Y)$  this Fourier-Mukai type transform. Then $F$ has a right adjoint $G$. \\
\end{Lem}

\textbf{Proof}. 
 We use the following isomorphisms: 
{\footnotesize
\begin{eqnarray} 
\Hom_{D_{qc}(Y)}( \R         p_{2\;*}\cL p_{1}^{*}(A) \bL \mathcal{E}, B) &\cong& 
\Hom_{D_{qc}(X \times Y)}(\cL p_{1}^{*}(A) \bL \mathcal{E}, p_{2}^{!}B) \label{4.1}\\
&\cong& 
\Hom_{D_{qc}(X \times Y)}( \cL p_{1}^{*}(A), \R        \mHom(\mathcal{E}, p^{!}_{2}B)) \label{4.2} \\
&\cong& 
\Hom_{D_{qc}(X)}( A, \R p_{1\;*} \; \R \mHom(\mathcal{E}, p^{!}_{2}B)). 
\end{eqnarray}}
The isomorphism \eqref{4.1} follows from  Grothendieck duality (see \cite{Ne96:gnus}). The isomorphism \eqref{4.2} follows from the  fact that  ($\otimes  $, $\; \mHom$) is an adjoint pair. The last isomorphism is a consequence of the fact that 
 ( $\cL p_{1} ^{*}$, $\; \R p_{1 \;*}$)  is an adjoint pair.
Thus $F$ has a right adjoint $G$. 
$\Box$
\begin{Remark}\label{boundbelow}
When the object $\mathcal{E}$ satisfies the assumptions in 
 Lemma~\ref{finite}, we have  $G(b) \in D^{+}(X)$ for all $b \in D^{b}(Y)$ by the explicit formula of the right adjoint $G$.
\end{Remark}

\subsection{Conclusion of the proof}
\textbf{Proof of Proposition~\ref{fin}}. 
By Lemma~\ref{rig}, the Fourier-Mukai type  transform $F:D_{qc}(W) \rightarrow D_{qc}(X)$ has a right adjoint $G: D_{qc}(X) \rightarrow D_{qc}(W)$,  
 so we   have the natural transforms 
$id_{D_{qc}(W)} \rightarrow GF$ and $FG \rightarrow id_{D_{qc}(X)}$. 
To show $F: D^{b}(W) \rightarrow D^{b}(X)$ 
is an equivalence, it suffices 
to show that $a \cong GF(a)$ for all
 $a \in D^{b}(W)$ and  $FG(b) \cong b$ for all $b \in D^{b}(X)$.


For each $a \in D^{b}(W)$ we have a distinguished  triangle in $D_{qc}(W)$ 
\[ \rightarrow  a \rightarrow GF( a) \rightarrow c \rightarrow a[1] \rightarrow \; . \; \; \;(*)  \]
To show that $a \cong GF(a)$,  it amounts  to showing  $c \cong 0$. We first show a weaker claim.

\begin{cl}\label{cbound}
$c \in D^{b}(W)$.

\end{cl}

Note that Claim~\ref{cbound} is equivalent to the fact that $GF(a) \in D^{b}(W)$. 
Pulling back everything to each $Y_{i}$, 
we get a  distinguished  triangle in $D_{qc}(W_{i})$    
\[\rightarrow  a_{i} \rightarrow G_{i}F_{i}(a_{i}) \rightarrow c_{i} \rightarrow a_{i}[1] \rightarrow \;  \; \; \; (*)_{i}\]
for each $Y_{i}$.   

Note that for every $x \in D^{b}(X)$ we have $G(x) \in D^{+}(W)$ by the explicit formula of the right adjoint functor, so $GF(c) \in D^{+}(W)$ (see Remark~\ref{boundbelow}). 

Since $F_{i}:D^{b}(W_{i}) \rightarrow D^{b}(X_{i})$ is 
an equivalence by assumption, 
it follows that $\Hom^{j}_{D_{qc}(W_{i})}(x_{i}, c_{i})=0$ for 
all $j$ and  all $x_{i} \in D^{b}(W_{i})$. In fact, we only need $F_{i}$ to be fully faithful for this assertion.  
To show $c_{i} \cong 0$, we use the following triangle for each $k$
\[ \rightarrow \tau_{\leq k}c_{i} \rightarrow c_{i} \rightarrow \tau_{\geq k+1}c_{i} \rightarrow \tau_{\leq k}c_{i}[1] \rightarrow \; .\]
Since $c_{i} \in D^{+}(W_{i})$, it follows that $\tau_{\leq k}c_{i} \in D^{b}(W_{i})$ for all $k$. Taking $\Hom(\tau_{\leq k}c_{i},-)$ into the above triangle, and
 noticing  that $\Hom^{0}_{D_{qc}(W_{i})}(\tau_{\leq k}c_{i},\tau_{\geq k+1}c_{i})=0$ and $\Hom^{-1}_{D_{qc}(W_{i})}(\tau_{\leq k}c_{i},\tau_{\geq k+1}c_{i})=0$, it follows that
 $\Hom^{0}_{D_{qc}(W_{i})}( \tau_{\leq k}c_{i}, \tau_{\leq k}c_{i}) \cong \Hom^{0}_{D_{qc}(W_{i})}( \tau_{\leq k}c_{i}, c_{i})$, which is $0$ since $\tau_{\leq k}c_{i} \in D^{b}(W_{i})$. 

If $c_{i} \not\cong 0$, then  we can choose a $k$ such that 
$\tau_{\leq k}c_{i} \not\cong 0$. For such a $k$, we have 
$\Hom^{0}_{D_{qc}(W_{i})}( \tau_{\leq k}c_{i}, \tau_{\leq k}c_{i}) \neq 0$, a contradiction.
This shows that $c_{i} \cong 0$. In particular, $c_{i} \in D^{b}(W_{i})$, so 
$c \in D^{b}(W)$. This proves Claim~\ref{cbound}. \newline

Let $\Omega$ be as in Lemma~\ref{spa}.   Let $y \in \Omega$. Taking $\Hom(y,-)$ into the distinguished triangle $(*)$ and the distinguished  triangle $(*)_{i}$ for each $Y_{i}$ , we get the following exact sequences 
\newline  {\scriptsize
\xymatrix{ \Hom^{j}(y,a)\ar[r] \ar[d] &  \Hom^{j}(y,GF(a)) \ar[r] \ar[d]& \Hom^{j}(y,c) \ar[r] \ar[d] &  \Hom^{j+1}(y,a) \ar[r] \ar[d] &  \\
    \Hom^{j}(y_{i},a_{i})\ar[r] &  \Hom^{j}(y_{i},G_{i}F_{i}(a_{i})) \ar[r] & \Hom^{j}(y_{i},c_{i}) \ar[r]  &  \Hom^{j+1}(y_{i},a_{i}) \ar[r]  & _{\;.} \\ } }\newline
\newline 
Note that the support of $y$ lies  in some $Y_{i}$ since $y \in \Omega$. Fix such a 
scheme $Y_{i}$. 
We have $y_{i} \cong y$  and all  vertical arrows are isomorphisms.  
Since $(F,G)$ and $(F_{i},G_{i})$ are adjoint pairs, it follows that
 $\Hom^{j}_{D_{qc}(W)}(y, GF(a)) \cong \Hom^{j}_{D_{qc}(X)}(F(y),F(a))$ and  $\Hom^{j}_{D_{qc}(W_{i})} (y_{i},G_{i}F_{i}(a_{i})) \cong \Hom^{j}_ {D_{qc}(X_{i})}(F_{i}(y_{i}),F_{i}(a_{i}))$. Together with our assumption  that all $F_{i}$ are equivalences, 
 this implies $\Hom^{j}_{D_{qc}(W)}(y,c)=0$ for all $y \in \Omega$. 
Since $c \in D^{b}(W)$ and $\Omega$ is a spanning class for $D^{b}(W)$, it 
follows that  $c \cong 0$. \newline

To show $FG(b) \cong b$ for all $b \in D^{b}(X)$,  we need to use the assumption that $F_{i}$ is an equivalence. 
Note that from Claim~\ref{cbound} we have $G_{i}F_{i}(a_{i}) \in D^{b}(W_{i})$ for all $a_{i} \in D^{b}(W_{i})$, and  since $F_{i}:D^{b}(W_{i}) \rightarrow D^{b}(X_{i})$ is an equivalence it follows that $G_{i}: D^{b}(X_{i}) \rightarrow D^{b}(W_{i})$.   
Therefore $(F_{i}, G_{i})$ is also an adjoint pair when we work on $D^{b}$,   so  $G_{i}:D^{b}(X_{i}) \rightarrow D^{b}(W_{i})$ is also an equivalence.
Using  another distinguished triangle 
\[ \rightarrow FG(b) \rightarrow b \rightarrow c \rightarrow FG(b)[1] \rightarrow ,\; \]
 one can show that  $FG(b) \cong b$ by a similar argument. 
This concludes the proof. $\Box$

\section{Basic properties of $M(X/Y)$ and $W(X/Y)$}\label{sec:bas}

We  prove some basic lemmas on  the distinguished component $W$ of $M(X/Y)$.   
\subsection{Characterization of the universal perverse ideal sheaf}
We begin with the next lemma.
\begin{Lem}\label{torsionfree}
Let $S$ and $T$ be two integral schemes and $F$ be a coherent  sheaf on $S \times T$. Let  $\pi$ be the projection map $S \times T \rightarrow T$. 
Assume the following two conditions: 
\begin{enumerate}[{(}1{)}]
\item $F$ is flat over $S$, and
\item there is a dense open set $U \subset S$ such that $F$ is torsion free on $\pi ^{-1}(U)$.

Then the sheaf $F$ is torsion free.

\end{enumerate}  

\end{Lem}
\textbf{Proof}.
The problem is local, so we may assume that both $S$ and $T$ are affine schemes. We  use torsion sections to get a contradiction.

Assume that $F$ is not torsion free. Let $x$ be a torsion section.  Denote by $V(y)$ the zero scheme of $y$ for a regular function $y$ on $S$.   By the assumption (2),  the image of the support of $x$ under the projection,  denoted  by  $\pi_{S}({\rm Supp}(x))$, is a proper subscheme of $S$.  We can find a regular function $s$ on $S$  such that $\pi_{S}({\rm Supp}(x)) \subset V(s)$ and the regular function $s$ annihilates $x$ (we consider $s$ as a regular function on $S \times T$ by the natural map of rings induced by the projection map). 
 
Consider the exact sequence \newline
\xymatrix{ 0 \ar[r]& \mathcal{O}_{S \times T} \ar[r]^{\cdot \;s} & \mathcal{O}_{S \times T} \ar[r] &\mathcal{O}_{V(s) \times T} \ar[r]& 0.}\newline
Tensoring this with $F$, we get a right exact sequence \newline
\xymatrix{& F \ar[r] & F \ar[r]& F \otimes \mathcal{O}_{V(s) \times T} \ar[r]& 0.}\newline
The map on the left is the multiplication by $s$. Since $xs=0$, it is not injective. This shows that $Tor_{1}(F, \mathcal{O}_{V(s) \times T}) \neq 0$. This implies that $F$ is not flat, a contradiction. $\Box$

We give a proposition on the universal ideal sheaf.

\begin{Prop}\label{univ}
The universal perverse ideal sheaf is  the 
ideal sheaf $\mathcal{I}_{W \times _{Y} X}$ of the fiber product, consequently the universal perverse point sheaf is $\mO_{W \times _{Y} X}$.

\end{Prop}
\textbf{Proof}.
Let $F$ be the universal ideal sheaf 
and $\alpha:F \rightarrow \mathcal{O}_{W \times X}$ 
be the corresponding  homomorphism between sheaves. 
Denote by $\Gamma$ 
the graph of $g: W \rightarrow Y$. 
The sheaf $F$ is flat over $W$ by definition. 
It is clear that $F$ is torsion free on the dense open set $U \times X$, where $U$ is the isomorphic locus of $f: X \rightarrow Y$ and is considered as an open set inside both $X$, $Y$ and $W.$
 
By Lemma~\ref{torsionfree}, it follows that $F$ is indeed torsion free. Since the morphism 
$\alpha : F \rightarrow \mathcal{O}_{W \times X}$ 
is generically injective,   the kernel is a torsion subsheaf. By Lemma~\ref{torsionfree} again, it follows that the homomorphism $\alpha$ is injective.
 So we can identify $F$ as an ideal sheaf of $\mathcal{O}_{W \times X}$. 

We show that $F= \mathcal{I}_{W \times _{Y} X}$. 
As shown in  \cite{br:gnus}, 
we have that $f_{*}(F)= I_{\Gamma}$, the ideal sheaf of the graph in $W \times Y$.   
By Proposition 5.1 in \cite{br:gnus} , 
it follows that the natural 
map $f^{*}f_{*}(F) \rightarrow F$ is surjective. 

Since  $f_{*}(F)= I_{\Gamma}$, 
the images of 
$f^{*}f_{*}(F)$ and $f^{*}f_{*}(\mathcal{I}_{W \times _{Y} X})$ 
in  $\mathcal{O}_{W \times X}$ coincide. 
This shows that $F = \mathcal{I}_{W \times _{Y} X}$. $\Box$
\subsection{Flatness lemma}
\begin{Prop}\label{flat}
Let $X_{1}$ be an irreducible quasi-projective  variety.   Consider the diagram : \newline
\xymatrix{ &&&X_{1} \times _{Y} X \ar[r] \ar[d] & X \ar[d] \\
           &&&X_{1} \ar[r]_{f_{1}} & Y . \\}\newline   
If the ideal $\mathcal{I}_{ X_{1} \times _{Y} X}$ in $\mathcal{O}_{X_{1} \times X}$ is flat over $\mathcal{O}_{X_{1}}$ and the image of $X_{1}$ is not contained in the image of the exceptional set of $X/Y$ in $Y$, then there is a canonical morphism $h:X_{1} \rightarrow W(X/Y)$. 
\end{Prop}

\textbf{Proof}.
Let $U$ be the isomorphic locus of $X \rightarrow Y.$ We consider $U$ as an open subset both in $X$ and $Y$.
Pick a point $x_{1} \in X_{1}$ such that $u=f_{1}(x_{1}) \in U$.
The sheaf $\mathcal{I}_{X_{1} \times _{Y} X , \; x_{1}}$ is 
$\mathcal{I}_{u}$, the ideal of the point  $u \in U.$ 
Since $\mathcal{I}_{X_{1} \times _{Y} X}$  is flat, this family of sheaves
 has the correct numerical class, say $\gamma .$
 
The scheme $W(X/Y)$ is isomorphic to $M_{PI}(X/Y, \gamma)$, the moduli space of perverse ideal sheaves   with the numerical equivalence class $\gamma$.  
It suffices to show that there is a morphism $h: X_{1} \rightarrow M_{PI}(X/Y, \gamma)$, which  amounts  to showing  that $\mathcal{I}_{X_{1} \times _{Y} X}$ is a family of perverse ideal sheaves.  This would follow if we can show 
that  the natural homomorphisms  
\[ f^{*}f_{*}(\mathcal{I}_{X_{1} \times _{Y} X,\;x_{1}}) \rightarrow \mathcal{I}_{X_{1} \times _{Y} X,\;x_{1}} \] 
are surjections for all $x_{1} \in X_{1}$.
 This holds
if the natural homomorphism 
 \[f_{X_{1}}^{*}f_{X_{1} \;*}(\mathcal{I}_{X_{1} \times _{Y} X}) 
 \rightarrow \mathcal{I}_{X_{1} \times _{Y} X} \] 
is a  surjection, which follows  since $f_{X_{1} \;*}( \mathcal{I}_{X_{1} \times _{Y} X})$ is the ideal of the graph $f_{1}: X_{1} \rightarrow Y$. $\Box$

\subsection {Relations between $W(X/Y)_{T}$ and $W(S/T)$}
Let $X \rightarrow Y$ be a flopping contraction between three-dimensional normal varieties. Assume that $X$ has  at worst terminal Gorenstein singularities. By standard results on flops, the variety $Y$ has at worst terminal Gorenstein singularities (see Theorem 6.14 in \cite{km:gnus}). 
Let   $T \subset Y$ be an effective Cartier divisor; for simplicity we assume  that it  is an integral subscheme of $Y.$ 
Consider the diagram \newline
\xymatrix {&&& S \ar[d]_{f_{T}} \ar[r]_{i_{S}} & X \ar[d]_f \\
            &&&T            \ar[r]_{i_{T}}            & Y.}\newline    
Note that the conditions (B.1) and (B.2) hold for the morphism 
$S \rightarrow T$. 

The condition (B.2) is clear.  
We now show the condition (B.1).
It is clear that $f_{T \;*}(\mO_{S})= \mO_{T}$.  To show $\R ^{i}f_{T}(\mO_{S})=0$ for all $i \geq 1$, we apply the theorem on formal functions (see p.277 in \cite{ha77}). 
It suffices to show that 
\[   H^{i}(S_{t}, \mO_{t})=0 \; \; \forall \; i  \geq 1 \] 
for all $t \in T$.
 Since 
\[ \R ^{i}f(\mO_{X})=0 \; \; \forall \; i \geq 1 \]
by assumption, it follows that 
\[   H^{i}(X_{y}, \mO_{y})=0 \; \; \forall \; i \geq 1 \] 
for all $y \in Y.$
For any $t \in T \subset Y$ the fibers $X_{t}$ and $S_{t}$ are canonically isomorphic to each other. This implies that
\[ \R ^{i}f_{T}(\mO_{S})=0 \; \;  \forall \; i \geq 1. \] 
Thus the condition (B.1) also holds for $S \rightarrow T.$ 

\begin{Prop}{\rm (= Proposition~\ref{goingdown1})}\label{dim3}
There is a natural embedding $ W(X/Y)_{T} \hookrightarrow M(S/T)$, which is an inclusion of components. 

\end{Prop}

\textbf{Proof}. (a) We show that there is a canonical morphism $M(S/T) \rightarrow M(X/Y)_{T}.$ \newline

Let $p \in M(S/T).$ Denote the corresponding perverse point sheaf for $S \rightarrow T$ by $E_{p}.$    
It is clear that if for a point  $p \in M(S/T)$, the corresponding object $E_{p}$ is also a  perverse point sheaf for $X \rightarrow Y$, then this point, which
 we still denote by $p$, must lie in the fiber $M(X/Y)_{T}.$ 
Let
\[ 0 \rightarrow I_{E_{p}} \rightarrow \mO_{S} \rightarrow E_{p} \rightarrow 0 \]
be the exact sequence in the abelian category $Per(S/T).$ \newline
\newline $\textbf{Step 1}\; $ We show that for every point $p \in M(S/T)$, the corresponding perverse point sheaf $E_{p}$ for $S \rightarrow T$ is a perverse sheaf for $X \rightarrow Y.$ \newline

This follows easily by checking the conditions (PS.1)-(PS.3) of Lemma~\ref{ps1-3}.\newline
\newline $\textbf{Step 2}\; $ We show that $I_{E_{p}}$ is also a perverse sheaf. \newline

This again follows by checking the conditions (PS.1)-(PS.3). \newline 

Combining results from Step 1 and  Step 2, it follows that
\[ 0 \rightarrow I_{E_{p}} \rightarrow \mO_{S} \rightarrow E_{p} \rightarrow 0 \]
is also an exact sequence in the abelian category $Per(X/Y).$ \newline
\newline $\textbf{Step 3}\;$ The sheaf $\mO_{S}$ is a perverse structure sheaf for $X \rightarrow Y.$ \newline

Consider the exact sequence of sheaves
\[ 0 \rightarrow I_{S} \rightarrow \mO_{X} 
\rightarrow \mO_{S} \rightarrow 0. \]

It suffices to check that $I_{S}$ is a perverse ideal sheaf. This follows since  the conditions (PIS.1) and (PIS.2) of Proposition~\ref{pi1-2} are satisfied. 

Composing   two surjections $\mO_{X} \rightarrow \mO_{S}$ and 
$ \mO_{S} \rightarrow E_{p}$, we obtain the surjection  $\mO_{X} \rightarrow E_{p}$ in the abelian category $Per(X/Y).$ This shows that $E_{p}$ is a perverse point sheaf for $X \rightarrow Y.$ 
Since $M(X/Y)$ is a fine moduli space, we have an embedding  $M(S/T) \rightarrow M(X/Y)_{T}.$ \newline

(b) We show that there is an embedding $W(X/Y)_{T} \rightarrow M(S/T).$ \newline

For each point $w \in W(X/Y)_{T}$, let $E_{w}$ be  the corresponding perverse point sheaf. \newline
\newline $\textbf{Step 1}\;$ We  prove that $E_{w}$ is a perverse sheaf for $S \rightarrow T.$ \newline 

 The main point is that $E_{w}$ is indeed a complex of 
$\mO_{S}$-modules since the universal perverse point sheaf is the structure sheaf of the fiber product $W \times _{Y} X$. 
By checking the conditions (PS.1)-(PS.3)  of  Lemma~\ref{ps1-3} in Appendix~\ref{sec:per}, it follows that $E_{w} \in Per(S/T)$. \newline
\newline $\textbf{Step 2}\;$ The sheaf $\mO_{S}$ is a perverse structure sheaf for $X \rightarrow Y.$ \newline

This is proven in part (a). \newline
\newline $\textbf{Step 3}\;$ The sheaf $E_{w}$ is a perverse point sheaf for $S \rightarrow T$. \newline

The morphism $\mO_{X} \rightarrow E_{w}$ factors through $\mO_{S}$. By 
Step 2 the morphism  $\mO_{X} \rightarrow \mO_{S}$ is a surjection in the abelian category $Per(X/Y)$, so  $\mO_{S} \rightarrow E_{w}$ is also a surjection by standard results on abelian categories, which  shows  the  corresponding kernel,  denoted  by $I_{E_{w}}$, is also a perverse sheaf. Note that the object $I_{E_{w}}$
 is also a shifting of the cone of $\mO_{S} \rightarrow E_{w}$, from which 
 follows that  $I_{E_{w}}$ is  a complex of $\mO_{S}$-modules. 
Abusing the notation,
 we  denote by $I_{E_{w}}$ and $E_{w}$ the objects  in
 $Per(S/T)$  such that
 $I_{E_{w}} \rightarrow \mO_{S} \rightarrow E_{w}$ is a
 distinguished triangle in $D^{b}(S)$ and the push-forward of this triangle is the exact sequence 
\[ 0 \rightarrow I_{E_{w}} \rightarrow \mO_{S} \rightarrow  E_{w} \rightarrow 0 \]  
in  the abelian category $Per(X/Y)$. 
Since $E_{w}$ is in the correct numerical class,
it follows that $E_{w}$ is a perverse point sheaf for $S \rightarrow T$.
 This gives an embedding $W(X/Y)_{T} \hookrightarrow W(S/T)$ since $W(S/T)$
  is a fine moduli space. \newline

Combining the results in (a) and (b), we obtain two morphisms 
$M(S/T) \rightarrow M(X/Y)_{T},$ and $W(X/Y)_{T} \rightarrow M(S/T).$
Each of these two morphisms is an embedding. By our  construction, the composition $W(X/Y)_{T} \rightarrow M(S/T) \rightarrow M(X/Y)_{T}$ is an inclusion of components, which implies  that each morphism is an inclusion of components.
This concludes the proof. $\Box$ \newline

Since $W(X/Y)$ is a fine moduli space, pulling back the universal object over $W(X/Y)$ via the canonical embedding $W(S/T) \rightarrow W(X/Y)$ in part (a), one obtains the following corollary of Proposition~\ref{dim3}.

\begin{Cor}\label{universalobject}
The universal perverse point sheaf 
for the morphism $S \rightarrow T$  is 
 $ \mO_{W(S/T) \times_{T} S} \cong    \cL  i^{*} ( \mO_{W(X/Y) \times_{Y} X})$.   $ \Box$ 
\end{Cor} 

\section{Deformations and general hyperplane sections}\label{sec:defo}

The proof in this section was inspired by helpful discussions with M. Reid. 
 Throughout  this section we assume that $Y$ is an  affine variety. 
Let $f:X \rightarrow Y$
be a crepant projective birational  morphism between two quasi-projective three-dimensional normal Gorenstein varieties. Assume that the variety $X$ has at worst terminal  singularities.  Denote  the exceptional set by $C.$  
Under these assumptions, all  singularities of $X$ are isolated  hypersurface singularities (see \cite{km:gnus} p.169). 
By standard results in the MMP, it is well-known  that $Y$ is also terminal (see Theorem 6.14 in \cite{km:gnus}).  

We first show  that for  a general one-parameter deformation $F: \mathcal{X} \rightarrow \mathcal{Y}$ of $f:X \rightarrow Y$ the total space $\mathcal{X}$ is nonsingular. Then we show that the hyperplane section $\mS$,  the preimage  of a general member $\mathcal{T}$ of a suitable linear system of divisors passing through the singular points $Y_{sing}=  \{ p_{i} : i=1, \cdots , m \} \subset Y$, is nonsingular. 
In the first part  we  use  the fact 
that these singularities are hypersurface singularities. 
The second part can be reduced to showing that the preimage of a general hyperplane passing through $Y_{sing} \in Y$ has only canonical singularities.   

Let $V_{0} \subset H^{0}(Y, \mathcal{O}_{Y})= H^{0}(X, \mathcal{O}_{X})$ be 
 any  linear sub-system of divisors passing through $Y_{sing}= \{ p_{i} : i= 1, \cdots ,m \}$ such that $|Bs V_{0}| = Y_{sing}$ (as a scheme). Let $T$ be  a general element of  $V_{0}$. Denote the preimage of $T$ in $X$  by $S$. 
\begin{Prop}\label{canonical}
The preimage $S$ of a general element  $T$ of  $V_{0} \subset H^{0}(Y, \mathcal{O}_{Y})=H^{0}(X, \mathcal{O}_{X}) $  has only canonical singularities.
\end{Prop}

\textbf{Proof}.
 First note that  $S$ is a Gorenstein variety since it is a hyperplane section of  a Gorenstein variety $X$. The divisor $T$ has only canonical singularities.
By a Bertini type theorem, the Cartier divisor  $T$ is nonsingular outside  $Y_{sing}.$ 
Therefore $S$ is  nonsingular outside the 
exceptional curves $C$. 

We  show that $S$ is normal and has only canonical (Du Val) singularities. 
We have $K_{S}= f^{*}K_{T}$
by 
\begin{enumerate}[{(}1{)}]
\item $K_{S}=K_{X}|_{S} +S|_{S}$ and $K_{T}=K_{Y}|_{T}+T|_{T}$  (by the adjunction formula),  and 

\item $K_{X}=f^{*}K_{Y}$.

\end{enumerate}
Consider the normalization $g: S^{n} \rightarrow S$.
 We have  $\omega_{S^{n}}=( \mathcal{C}) g^{*}(\omega_{S})$, where $\mC$ is the conductor ideal. 
Since $T$ has only Du Val singularities, we have $(g \circ f)^{*} (\omega_{T})
\subset \omega_{S^{n}}$. This shows that $S$ is normal. 
To complete the proof, we  compute the discrepancies. Take a resolution $h: V \rightarrow S$ of $S$. We have 
\[ K_{V} = h^{*}K_{S} + \sum a_{i}E_{i} = (h \circ f)^{*}K_{T} + \sum a_{i}E_{i}  \]
where $E_{i}$'s are the exceptional divisors. 
Since $T$ has only  canonical singularities, it follows that $S$  has only 
 canonical singularities.  $\Box$

\begin{Prop}{\rm (= Proposition~\ref{smo1})}\label{smo}
A general one-parameter deformation   of $f: X \rightarrow Y$
is nonsingular.
\end{Prop}

\textbf{Proof}.
 Let $\mathcal{X}_{univ}$ be  the semiuniversal object over the semiuniversal deformation space $Def(X)$. Let $\mathcal{Y}= Spec ( \mathcal{O}_{\mathcal{X}_{univ}})$. Then $ \mathcal{Y}$ is a deformation of $Y$, and  hence the natural morphism $F: \mathcal{X}_{univ} \rightarrow \mathcal{Y}$ is a deformation of $f: X \rightarrow Y.$ 
 Thus it suffices   to deform  a Zariski neighborhood of $f^{-1}(p)$ in $X$. 

Since $X$ has only isolated hypersurface singularities, the deformation space of $X$ is $\Ext^{1}( \Omega_{X}, \mathcal{O}_{X})$. We show below that the obstruction group $\Ext^{2}( \Omega_{X}, \mathcal{O}_{X}) =0$. 
To compute $\Ext^{2}( \Omega_{X}, \mathcal{O}_{X})$, we use the following spectral sequence 
\[ H^{p}(X, \mExt^{q}(\Omega_{X}, \mathcal{O}_{X})) \Rightarrow \Ext^{p+q}(  \Omega_{X}, \mathcal{O}_{X}). \]
Since $X$ has only isolated hypersurface singularities, it is clear that $\mExt^{2}(\Omega_{X}, \mO_{X})$ is $0$. We also know that $H^{1}(X, \mExt^{1}(\Omega_{X}, \mO_{X}))=0$ since ${\rm Supp}(\mExt^{1}(\Omega_{X}, \mO_{X}))$ is isolated. It remains  to show that 
$H^{2}(X, \mExt^{0}( \Omega_{X}, \mathcal{O}_{X}))=0.$ 

This follows from the Leray spectral sequence 
\[ H^{p}(Y, \R  ^{q}f_{*}(\mathcal{F})) \Rightarrow H^{p+q}(X, \mathcal{F}) \]and $H^{i}(Spec(A), \mathcal{F})=0$ for $i \geq 1$. 
By a similar argument, one could obtain   that  $E^{p,q}_{2} =0$ for $p+q \geq 2$, though  we do not need this more general fact in our proof.

Since every $E^{p,q}_{2} =0$ for $p+q=2$, we get the following short 
exact sequence
\[ 0 \rightarrow H^{1}(X, \mHom( \Omega_{X}, \mathcal{O}_{X})) \rightarrow \Ext^{1}(\Omega_{X}, \mathcal{O}_{X}) \rightarrow H^{0}(X, \mExt^{1}(\Omega_{X}, \mathcal{O}_{X})) \rightarrow 0. \]

The important point is that the map $\Ext^{1}(\Omega_{X}, \mathcal{O}_{X}) \rightarrow H^{0}(X, \mExt^{1}(\Omega_{X}, \mathcal{O}_{X}))$ is surjective. Thus every deformation of the singularity can be lifted to a deformation of $X$.  Since $X$ has finitely many singularities and that smoothness at a given point is an open condition, it suffices to check the smoothness statement in  neighborhoods of each   singular point of $X$.    

Note that we can check whether a variety is  nonsingular at a given point $x$  locally analytically. Thus we shall work locally analytically  in the remainder of this argument.
Denote the semiuniversal  deformation space of the singularity  $x \in X$  by $Def(x \in X)$ and the  semiuniversal object over $Def(x \in X)$ by $\mathcal{X}$.
For an isolated hypersurface singularity, the total space $\mathcal{X}$  over the semiuniversal deformation space $Def(x \in X)$ is nonsingular by the explicit description of the semiuniversal space and the total space.

The variety $\mathcal{X}$ is analytically isomorphic to $f(x,y,z,w)+t_{1}f_{1}+ \cdots +t_{n}f_{n}=0$ where $n$ is the dimension of $Def(x \in X)$ and $f_{i}$ are suitable polynomials such that at least one of  the $f_{i}$, say $f_{1}$, is nonzero at $(0,0,0,0)$.

The canonical  morphism $\{0 \in Def(X) \} \rightarrow \{ 0 \in Def(x \in X)\}$ gives a linear map on tangent spaces. This map is surjective.   We write the defining equation for the semiuniversal object $\mathcal{X}_{univ}$  as $F= f(x,y,z,w)+ t_{1}g_{1}  + \cdots +t_{m}g_{m}=0$. Choose a direction of $c=(c_{1},\cdots,c_{m})$ in the tangent space of $Def(X)$ such that its image under the induced linear map is $(1,0,\cdots,0)$. This gives a smoothing of the singular point $x$ by the  Jacobi criterion. $\Box$ \newline

Let $F:\mX \rightarrow \mY$ be a one-parameter  deformation of $f:X \rightarrow Y$ such that $\mX$ is smooth. Let $\mathcal{T}$ be a general hyperplane section passing through $Y_{sing} \in Y$. Denote by  $\mathcal{S}$ the preimage of $\mathcal{T}$. 

\begin{Prop}\label{gen}
There exists a finite-dimensional vector   space $V \subset H^{0}(\mathcal{Y}, \mathcal{O}_{\mathcal{Y}})$  such that the preimage $\mathcal{S}$  of
  a general hyperplane section passing through $p$   is nonsingular. 
\end{Prop}
\textbf{Proof}.
Fix any  
 $V_{0}$ satisfying the conditions at the beginning of this section.
Let $V$ be  a finite-dimensional linear sub-system $V \subset H^{0}(\mathcal{Y}, \mathcal{O}_{\mathcal{Y}}).$ Denote by $Im(V)$ the image of $V$ in 
 $H^{0}(Y, \mO_{Y})$ under the natural homomorphism 
$H^{0}(\mathcal{Y}, \mathcal{O}_{\mathcal{Y}}) \rightarrow H^{0}(Y, \mO_{Y}).$ 
We  choose a finite-dimensional linear sub-system $V \subset H^{0}(\mathcal{Y}, \mathcal{O}_{\mathcal{Y}})$   such that  $ Im(V) \subset H^{0}(Y, \mO_{Y})$ contains $V_{0}$  and   $|Bs V| \cap Y = Y_{sing}=\{p_{i}: i=1, \cdots , m \}.$ 
A general hyperplane section $S$ of $V_{0}$ has only canonical surface singularities  and hence  has only hypersurface singularities. The subset of this linear system $V$ such that the corresponding members are nonsingular at a specific point is an open set. There are only  finitely many singular points on $S$. 
Combining these two facts, it suffices 
to check  the corresponding open set is nonempty for each singular point. We divide the singular points of $S$ into two types. 

A point $x \in S_{sing}$ is called of type 1 if 
$ x \in S_{sing} \bigcap X_{sing}.$  A point 
$y \in S_{sing}$ is called of type 2 if 
$y \in S_{sing}/X_{sing}.$ 

We now  show that every section 
$s \in H^{0}(X, \mO_{X})$ can be 
lifted to a section in 
$H^{0}(\mX, \mO_{\mX})$. 
This follows easily from the exact sequence 
\[ {\footnotesize 0 \rightarrow  H^{0}(\mathcal{X},\mathcal{O}(-X)) \rightarrow H^{0}(\mathcal{X}, \mathcal{O}_{\mathcal{X}}) \rightarrow H^{0}(X, \mathcal{O}_{X}) \rightarrow  H^{1}(\mathcal{X},\mathcal{O}(-X)) } \]  
and the fact that  $H^{1}(\mathcal{X}, \mO_{\mX}(-X))=0$ (by the Leray spectral sequence and the fact that $\mathcal{Y}$ has only rational singularities).  
We still denote a lifting of $s$ by $s$. 

The variety $S \subset \mathcal{X}$ is a complete intersection. Denote  the ideal by $\mathcal{I}_{S}=(s, g) \subset \mathcal{O}_{\mathcal{X}}$. 

For a singular point $x$ of type 1,
we show that the divisor defined by $g$ is nonsingular near $x$. We prove this by computing the embedding dimension of $X$ at $x$. 
Passing to a formal or analytic neighborhood of $x \in \mathcal{X}$,  we may  assume that the ring of this formal neighborhood is $k[[x,y,z,w]]$. 
We have  $m_{x,\;S}/m^{2}_{x,\;S}=(x,y,z,w)/(m^{2}_{x,\;S},s,g)$. This vector space is of dimension $3$ since $S$ has a  canonical surface singularity at $x$. 
Since $x$ is a singular point of $X=\{s=0\}$, it follows that $s \subset m^{2}_{x,\;S}$, which implies $\{g=0\}$ is nonsingular at $x$. 

For a singular point  $y$ of type 2 in $S$,  the defining equation $s$ of $X$ is nonsingular at $y$. For a small enough $\epsilon$ the hyperplane section defined by $\epsilon \cdot g+s$ gives a divisor, which is nonsingular at $y$. $\Box$  

\section{Equivalences of derived categories: dimension $4$ to dimension $3$}\label{sec:going} 
The proof in this section is based on suggestions of T. Bridgeland. 
We again assume that $Y$ is an affine quasi-projective variety throughout this section.
Let $X$ be a quasi-projective threefold with only terminal Gorenstein singularities and $f: X \rightarrow Y$ be a flopping contraction.  Let  $W=W(X/Y)$ be the distinguished component of the moduli space of perverse point sheaves $M(X/Y)$.  We  prove in Section~\ref{sec:defo} that there is  a deformation $F: \mX \rightarrow \mY$ of $f: X \rightarrow Y$ with smooth $\mX$. 
We summarize what we  know: 
\begin{enumerate}
\item The Fourier-Mukai type transform $\Psi : D^{b}(\mW) \rightarrow D^{b}(\mX)$  defined by the universal perverse point sheaf is an equivalence,  
\item  $\mW$ is smooth and is  the flop of $\mX \rightarrow \mY$, and   
\item  $W \cong \mW _{Y}$.    
\end{enumerate}
From what we know, it is a standard argument to  deduce  that $W(X/Y) \rightarrow Y$ is the flop of 
 $X \rightarrow Y$.
We sketch the argument here for the reader's convenience. 

Since $W$ is  Gorenstein and generically reduced, it is a reduced scheme and hence  is  an integral scheme. Using the argument given in Proposition~\ref{canonical}, it follows that   $W$ is normal and has at worst  terminal singularities.    
By the adjunction formula, 
we have $K_{W} \cdot \mC = K_{\mW} \cdot \mC$ for every curve $\mC \subset W \subset \mW$, which  implies the canonical bundle $K_{W}$ is $g$-trivial for $g: W \rightarrow Y$ since
 $K_{\mW}$ is $G$-trivial. Let $\mathcal{D}_{1} \subset \mX$ be the effective divisor such that $- \mathcal{D}_{1}$ is  $F$-ample and its birational transform $\mathcal{D}_{2}$ in $\mW$ is  
 $G$-ample.  Intersect  $\mathcal{D}_{1}$ with $X$, and  denote the intersection  by $D_{1}$. Then $-D_{1}$  is an $f$-ample divisor and $D_{2}$ is a $g$-ample divisor. To show $W$ is the flop, it remains to show  that the morphism $g$ is not a divisorial contraction, which is evident since $K_{W}=g^{*}K_{Y}$ and $Y$ has only terminal singularities. \newline

Our goal in this section  is to prove that
$\Psi_{0}: D^{b}(W) \rightarrow D^{b}(X)$ is an equivalence of categories (see below for the notation $\Psi_{0}$). 

Consider the  diagram \newline 
\xymatrix{ & W \times X \ar[dl]_{(p_{1})^{0}} \ar[dr]^{(p_{2})^{0}}&&&& \mathcal{W} \times_{C} \mathcal{X} \ar[dl]_{p_{1}} \ar[dr]^{p_{2}}&\\ 
          W \ar[dr]_{f_{0}}& & X \ar[dl]^{g_{0}} && \mathcal{W}\ar[dr]_f && \mathcal{X} \ar[dl]^g \\
         & Y\ar[d] &&&& \mathcal{Y}\ar[d] & \\
& \{0\}&& \ar[r]_{i_{0}} && \mathcal{C} \; . & }\newline
Denote by $\Psi_{0}:D^{b}(W) \rightarrow D^{b}(X)$ the Fourier-Mukai type transform defined by the kernel $\cL i^{*}_{0} ( \mO_{\mW \times _{\mY} \mX})$, which is equivalent to the Fourier-Mukai type transform defined 
by  the universal perverse point sheaf $\mO_{W \times_{Y} X}$ for $X \rightarrow Y$ (see Corollary~\ref{universalobject}). 

We claim  that $ \Psi ( i_{0\;*}  \mathcal{F}) \cong i_{0\;*} \circ \Psi_{0}(\mathcal{F})$ for $\mathcal{F} \in D^{b}(W).$ In fact, we  prove a stronger lemma below. \newline

Let $\mathcal{E} \in D^{b} (\mathcal{W} \times _{\mathcal{C}} \mathcal{X})$ be an object satisfying the assumptions in Lemma~\ref{finite}.
We may also consider it as an object in $D^{b}(\mW \times \mX)$.
Let $F: D^{b}(\mW) \rightarrow D^{b}(\mX)$ be the Fourier-Mukai type transform defined by $\mathcal{E}$.
By Lemma~\ref{fam}, the functor $F$ can be defined as $ \R p_{1\;*}(\cL p^{*}_{2}( (-) \bL \mathcal{E}))$.
Denote  by $F_{0}: D^{b}(W) \rightarrow D^{b}(X)$ the Fourier-Mukai transform defined by the object $\cL i_{0}^{*} \mathcal{E} \in D^{b}(W \times X)$. 
\begin{Lem}{\rm(= Proposition~\ref{gen1})}\label{commu}
Notation as above. Denote by $F$ the Fourier-Mukai type transform defined by the object $\mathcal{E}$. Then $F(i_{0\;*}(-)) \cong i_{0\;*} \circ F_{0}(-).$ 
\end{Lem}
\textbf{Proof}. 
We use the following isomorphisms: 
\begin{eqnarray}
F(i_{0\;*}(-)) &=&\R        p_{1\;*}(\cL p^{*}_{2}(i_{0\;*}(-) \bL  
                             \mathcal{E})) \nonumber  \\
               &\cong&    \R        p_{2\;*}(\R        i_{0\;*}(\cL (p_{1})^{0\;*}(-)) \bL \mathcal{E}) \label{6.1}\\
               &\cong&  
\R        p_{2\;*}(\R        i_{0\;*}(\cL (p_{1})^{0\;*}(-) \bL \cL i_{0}^{*}\mathcal{E})) \label{6.2} \\
&\cong& \R        i_{0\;*}(\R        (p_{2\;*})^{0}(\cL (p_{1})^{0\;*}(-) \bL \cL i_{0}^{*}\mathcal{E})). \label{6.3}
\end{eqnarray}
The  isomorphism \eqref{6.1} follows from the flat base change theorem. The isomorphism \eqref{6.2} follows from the projection formula. The isomorphism \eqref{6.3} is obvious. 
The last line is, by definition, the functor $i_{0\;*} \circ F_{0}(-)$. 
$\Box$

\begin{Prop}{\rm (= Proposition~\ref{fourto3})}\label{4to3}
Notation as above. Then
\[ \Psi: D^{b}(\mW) \cong D^{b}(\mX) \; \Longrightarrow \; \Psi_{0}: D^{b}(W) \cong D^{b}(X). \]

\end{Prop}
\textbf{Proof}.
 Applying Lemma~\ref{commu} to $\Psi$, it follows   that $\Psi (i_{0 \; *}(-)) \cong i_{0 \; *}(\Psi_{0}(-))$. 
 Let $\Phi:  D^{b}(\mX) \rightarrow D^{b}(\mW)$  
be the right adjoint functor of 
$\Psi: D^{b}(\mW) \rightarrow D^{b}(\mX)$, and
 $\mathcal{E}_{1}$ be the object corresponding to the Fourier-Mukai type transform $\Phi$. 

Denote by $\Phi_{0}: D^{b}(X) \rightarrow D^{b}(W)$ 
the Fourier-Mukai type  transform  
defined by the object $\cL i_{0}^{*}\mathcal{E}_{1}$.  
Lemma~\ref{commu}  also implies   
$\Phi (i_{0 \; *}(-)) \cong i_{0 \; *}(\Phi_{0}(-))$. 
These two facts give   the following commutative diagram \newline
\xymatrix{&&& D^{b}(W) \ar[r]_{i_{0}} \ar[d]_{\Psi_{0}} &  D^{b}(\mathcal{W}) \ar[d]_{\Psi} \\           
&&&D^{b}(X) \ar[r]_{i_{0}} \ar[d]_{\Phi_{0}} & D^{b}(\mathcal{X}) \ar[d]_{\Phi} \\
 &&&D^{b}(W)\ar[r]_{i_{0}}&  D^{b}(\mathcal{W}).  
} \newline  
Combining the top and the bottom parts  of this diagram, it follows that  
$\Phi \circ \Psi(i_{0 \;}(-)) \cong i_{0 \; *}(\Phi_{0} \circ \Psi_{0}(-))$. 
The functor $\Phi \circ \Psi$ is the Fourier-Mukai type transform defined by the diagonal $\Delta_{\mathcal{W}} \hookrightarrow \mathcal{W} \times_{\mathcal{C}} \mathcal{W}$ (see \cite{bkr:gnus} or Appendix~\ref{sec:2ndapp}), so it is equivalent to the identity functor  $id_{D^{b}(\mW)}$.  Since  $i_{0}$ is a  closed embedding,  $R^{i}i_{0 \;*}(-) =0$ for $i \neq 0$.  Therefore 
$\Phi_{0} \circ \Psi_{0}(\mathcal{F} ) \cong \mathcal{F}$ for all  objects  $\mathcal{F} \in D^{b}(W)$. \newline

To show $\Psi_{0}  \circ \Phi_{0} \cong id$, we first note that  $\Phi$ is  an equivalence when $\Psi$ is an equivalence.
By a similar argument, one can show
that $\Psi_{0}  \circ \Phi_{0} \cong id$. $\Box$

\appendix \section{Perverse coherent sheaves}\label{sec:per}
We give the definitions and related results of perverse coherent sheaves in this section. The main reference for this appendix is \cite{br:gnus}. 

Let $f: X \rightarrow Y$ be a projective birational   morphism between quasi-projective varieties. The following two assumptions  are the same as in \cite{br:gnus}:
\begin{enumerate}[{(B.}1{)}]\label{fiber}
\item $\R  f_{*} \mathcal{O}_{X} = \mathcal{O}_{Y}$, and 
\item   every  fiber of $f$ is of dimension at most $1$. 
\end{enumerate}
Any flopping contraction of a canonical threefold   satisfies these two conditions. 

We write $\mathcal{A} = D(X)$ and $\mathcal{B} = D(Y)$. By Proposition 2.3 in \cite{br:gnus}, we can identify $\mathcal{B}$ with a right admissible triangulated sub-category of $\mathcal{A}$. We thus have a semiorthogonal decomposition $( \mathcal{C} , \mathcal{B} )$ where  
\[\mathcal{C} = \mathcal{B}^{ \perp} = \{ E \in D(X): \R f_{*}(E)=0 \}.\]

\begin{Lem}
An object $E \in D(X)$
lies in $\mathcal{C}$ precisely when its cohomology sheaves $H^{i}(E)$ lie in $\mathcal{C}$.
\end{Lem}
\textbf{Proof}.
This is Lemma 3.1 in \cite{br:gnus}. The proof is an easy spectral sequence argument.
The condition (B.2)  is needed in the proof. $\Box$

Now we can get a $t-$structure on $\mathcal{A}$ by gluing the $t-$structures on $\mathcal{B}$ and $\mathcal{C}$ (see \cite{bb:gnus} 1.4.8-10).  
The standard $t-$structure on $\mathcal{A}$ induces a 
$t-$structure $\mathcal{C}^{\leq 0} = \mathcal{C} \cap \mathcal{A}^{\leq 0}$ on 
$\mathcal{C}$. Shifting this by $p$ and gluing it to the standard $t-$structure on $\mathcal{B}$ gives a new $t-$structure on 
$\mathcal{A}$. 

This $t-$structure has the following properties: 

$\mathcal{A}^{\leq 0}_{p} = \{ E \in   \mathcal{A} : \R f_{*}(E) \in \mathcal{B}^{\leq 0} \;and \; \Hom_{\mathcal{A}}(E,C)=0 \;for \;all \;C \in \mathcal{C}^{\geq p} \},$

$\mathcal{A}^{\geq 0}_{p} = \{ E \in   \mathcal{A} : \R f_{*}(E) \in \mathcal{B}^{\geq 0} \;and \; \Hom_{\mathcal{A}}(E,C)=0 \;for \;all \;C \in \mathcal{C}^{\leq p} \}.$
\\
The heart of this $t-$structure is an abelian category 
$Per^{p}(X/Y) = \mathcal{A}^{\leq 0}_{p} \cap \mathcal{A}^{\geq 0}_{p}$.\\
We shall only consider $p=-1$ and call this category $Per(X/Y)$. Following Bridgeland, the objects of $Per(X/Y)$ are called perverse coherent sheaves.

The next lemma gives an explicit description of $Per(X/Y)$.

\begin{Lem}\label{ps1-3}
An object $E$ of $D(X)$ is a perverse sheaf if and only if the following three conditions are satisfied: 
\begin{enumerate}[{\rm (PS.}\rm 1{)}]
\item $H_{i}(E)=0$ unless $i=0$ or $1$,
\item $\R^{1}f_{*}H_{0}(E)=0$ and $\R^{0}f_{*}H_{1}(E)=0$,
\item $\Hom_{X}(H_{0}(E),C)=0$ for any sheaf $C$ on $X$ satisfying $\R f_{*}(C)=0.$ 
\end{enumerate}
\end{Lem}
\textbf{Proof}. This is Lemma 3.2 in \cite{br:gnus}. $\Box$

\begin{Def}
Two objects $A_{1}$ and $A_{2}$ of $D^{b}(X)$ are $\textit{numerically equivalent}$ if for any locally-free sheaf $L$ on  $X$ we have $\chi(L,A_{1})= \chi(L,A_{2}).$

\end{Def}

\begin{Def}
An object $F$ of $D(X)$ is a perverse ideal sheaf if there is an injection $F \hookrightarrow \mathcal{O} _{X}$ in the abelian category $Per(X/Y)$. An object $E$ of $D(X)$ is a perverse structure sheaf if there is a surjection $\mathcal{O}_{X} \rightarrow E$ in  $Per(X/Y)$. A perverse point sheaf is a perverse structure sheaf which is numerically equivalent to the structure sheaf of a point $x \in X$.
\end{Def}
 
A perverse ideal sheaf $F$ determines and is determined by a perverse structure sheaf $E$, which fit in an exact sequence in $Per(X/Y)$

\xymatrix { &&0 \ar[r]& F \ar[r] & \mathcal{O}_{X} \ar[r] & E \ar[r]& 0. }  

It turns out that a perverse ideal sheaf is a sheaf.
We quote proposition 5.1 in \cite{br:gnus}.

\begin{Prop}\label{pi1-2}
A perverse ideal sheaf on $X$ is, in particular, a sheaf on $X$. A sheaf on $X$ is a perverse ideal sheaf if and only if the following two conditions are satisfied:

\begin{enumerate}[{\rm(PIS.}\rm 1{)}]
\item the sheaf $f_{*}(F)$ on $Y$ is an ideal sheaf, and
\item the natural map of sheaves $f^{*}f_{*}(F) \rightarrow F$ is surjective.
\end {enumerate}

\end{Prop}

 Let $S$ be a scheme. Given a point $s \in S$, let $j_{s} : {s} \times X \rightarrow S \times Y $ be the embedding. As indicated in Bridgeland \cite{br:gnus}, a family of sheaves on $X$ over $S$ can be characterized as  an object $\mathcal{F}$ of $D(S \times X)$ such that for each point $s \in S$ the object $\mathcal{F}_{s} = \cL         j_{s}^{*} ( \mathcal{F} )$ of $D(X)$ is a sheaf.

Following \cite{br:gnus}, we  define the moduli functor of perverse sheaves.  

\begin{Def}
A family of perverse sheaves on $X$ over a scheme $S$ is an object $\mathcal{E}$ of $D(S \times X)$ such that for each point $s \in S$ the object 
$\mathcal{E} _{s} = \cL         j_{s}^{*} ( \mathcal{F})$ of $D(X)$ is a perverse sheaf. Two such families $\mathcal{E}_{1}$ and $\mathcal{E}_{2}$ are equivalent if $\mathcal{E}_{2} = \mathcal{E}_{1} \otimes L$ for some line bundle pulled back from $S$. The moduli functor of perverse coherent sheaves assigns to each scheme $S$ the set of equivalence classes of perverse coherent sheaves on $X \times S$.  

\end{Def}

The following theorem can be found in [Br00]: 

\begin{Thm}
The functor which assigns to a scheme $S$ the set of equivalence classes of families of perverse point sheaves on $X$ over $S$ is  representable by a projective scheme $M(X/Y)$. 
\end{Thm}

The scheme $M(X/Y)$ has a distinguished irreducible component which is  birational to $Y$. We shall call it $W(X/Y)$. When no confusion is possible, we  denote it by $W$.

\begin{Remark}\label{quai}
In \cite{br:gnus} Bridgeland proved  
the existence of a fine  moduli space of perverse ideal sheaves  
when $X$ and $Y$ are projective varieties. 
We can generalize  his existence result   to quasi-projective varieties.   
A simple observation below  shows 
how to weaken the projectivity assumption on  $Y$.
 
Let $f: X \rightarrow Y$ be a projective morphism between 
two quasi-projective varieties 
satisfying the conditions (B.1) and (B.2). 
We can find a 
completion  of $f: X \rightarrow Y$ as in the following diagram \newline
\xymatrix{&&&& X \ar[r]_{i_{1}} \ar[d]_f & \overline{X} \ar[d]_{\overline{f}} \\
          &&&& Y \ar[r]_{i_{2}}          &  \overline{Y}.}
 
The problem is that such a 
compactification may no longer 
satisfy  conditions (B.1) and (B.2). 
But we can define another moduli 
functor which  parameterizes the pairs of 
a sheaf $F$ and a 
homomorphism $\alpha: F \rightarrow \mO_{\overline{X}}$  
satisfying 
conditions (PIS.1)-(PIS.2) in  Proposition~\ref{pi1-2} (see also Proposition 5.1 in \cite{br:gnus}), and  the condition that $\overline{f}_{*}(F)$ is the ideal of some point $y \in \overline{Y}$.

The proof of the  existence of such a  moduli space 
is the same as the proof of the existence of the moduli space 
of perverse point sheaves in \cite{br:gnus}. 
When  restricted  
to $Y$, this scheme  is  the moduli  
space of perverse point sheaves 
for $f: X \rightarrow Y$. It is also evident that this construction is local in $Y$.
\end{Remark}

\section{Proof of Theorem~\ref{bkr}}\label{sec:2ndapp}
 
We sketch the  
 proof of Theorem~\ref{bkr}  in this appendix. The argument  is taken from   [BKR99]. 
We adapt the notation from Theorem~\ref{bkr}. The statement that $W$ is the flop is an easy corollary of the result on the equivalence of derived categories (see \cite{br:gnus}). We shall omit its proof. 

Consider the diagram \newline
\xymatrix {
&&&&& W \times X \ar[dl]_{\pi_{W}} \ar[dr]^{\pi_{X}}  \\
&&&& W  && X.} \newline 
Let $\mathcal{P} \in D^{b}(W \times X)$ be the universal perverse point sheaf.
We define a functor (using results in \cite{sp88}) 
\[ \Psi = R \pi _{X \;*} ( \mathcal{P} \bL \pi _{W}^{*} (-)) : D_{qc}(W) \rightarrow D_{qc}(X). \]
It turns out that $\Psi$ sends  $D^{b}(W)$ to $D^{b}(X)$ 
(see Step 1 below). \newline
\newline $\textbf{Step 1}$ \newline

Each $\mathcal{P}_{w}$ has bounded homology sheaves. The variety  $X$ is nonsingular. These imply that $\mathcal{P}$ is of finite homological dimension. So we have $\Psi : D^{b}_{c}(W) \rightarrow D^{b}_{c}(X)$. \newline
\newline $\textbf{Step 2}$ \newline

We define another functor $\Upsilon : D_{qc}(X) \rightarrow D_{qc}(W)$ by 
\[ \Upsilon (-) = [ \R         \pi _{W_{*}} ( \mathcal{P}^{\vee} \otimes \pi ^{*}_{X} ( \omega _{X} )[n]) \bL \pi ^{*}_{X}(-)], \]
where $\mathcal{P}^{\vee}$ is the derived dual of $\mathcal{P}$.
Note that  this functor sends objects in $ D^{b}_{c}(X)$ to $D^{b}_{c}(W)$ since $ \mathcal{P}^{\vee}$ is of  finite homological dimension. 
We now restrict this functor to $D^{b}_{c}$. Then $\Upsilon:  D^{b}_{c}(X) \rightarrow D^{b}_{c}(W)$
 is left adjoint to $\Psi:D^{b}_{c}(W) \rightarrow D^{b}_{c}(X)$ as shown in [BKR99] (p.16). 
The composite functor $\Upsilon \circ \Psi$ is given by $\R         \pi _{2_{*}} ( \mathcal{Q} \bL \pi _{1}^{*} (-))$, where $\pi _{1}, \pi _{2} : W \times W \rightarrow W$  are the projections and $\mathcal{Q}$ is some object of $D^{b}_{c}(W \times W)$. 
 
If $i_{w} : w \times W \hookrightarrow W \times W$ be the embedding, then  $\cL  i_{w}^{*}( \mathcal{Q}) = \Upsilon \Psi \mathcal{O}_{w}$.
We have the following isomorphisms 
\[ \Hom^{i}_{D(W \times W)} ( \mathcal{Q} , \mathcal{O}_{w_{1},w_{2}}) = \Hom^{i}_{D(W)}( \Upsilon \Psi \mathcal{O}_{w_{1}}, \mathcal{O}_{w_{2}})\]
\[  = \Ext ^{i}_{X} ( \Psi \mathcal{O}_{w_{1}} , \Psi \mathcal{O}_{w_{2}} )= \Ext^{i}_{X} ( \mathcal{P}_{w_{1}}, \mathcal{P}_{w_{2}} ). \]
Each $\mathcal{P}_{w}$ is simple, so its support is connected and since $\R         f_{*} ( \mathcal{P}_{w} ) = \mathcal{O}_{y}$, where $y=g(w)$, it follows that $\mathcal{P}_{w}$ is supported on a fiber of $f$ over $y$. Since $f$ is crepant, we have $\mathcal{P}_{w} \otimes \omega  = \mathcal{P}_{w}$. For distinct $w_{1},\;  w_{2}$, Serre duality together with Lemma 3.6 in [Br00] implies that  $\Hom_{D_{c}}^{i} ( \mathcal{P}_{w_{1}}, \mathcal{P}_{w_{2}} )= 0$ unless $g( w_{1})= g( w_{2})$ and $ 1 \leq i \leq n-1$. Since $X \rightarrow Y$ is crepant, it follows that $\mathcal{P}_{w} \otimes \omega _{X} = \mathcal{P}_{w}$. \newline
\newline $\textbf{Step 3}$ \newline

We  prove in Step 2 that ${\rm h.d.}(\mathcal{Q}) \leq (n-1)-1=n-2$ when restricted to $W \times  W - \Delta_{W}$. 
We know that  ${\rm dim}\; (W \times _{Y} W) \leq n+1$ by assumption, and  ${\rm Supp}(\mathcal{Q})$ is contained in  $W \times _{Y} W$. Since we have ${\rm codim} (\mathcal{Q}) \geq n-1$, the intersection theorem implies   $\mathcal{Q} \cong 0$ outside the diagonal. 

Fix a point $w \in W$, put $E= \Upsilon \circ \Psi (\mathcal{O}_{w})$. We  prove above that $E$ is supported at the point $w$. 

\begin{cl}
$H_{0}(E)= \mathcal{O}_{w}$.
\end{cl}
The proof of this claim can be found in [BKR99] (p.18).    
 Corollary 5.3 in [BKR99]  then implies  that $E \cong \mO_{w}$ and $W$ is non-singular. 
Applying Theorem 2.3 in [BKR99],  it follows  that $\Psi: D^{b}_{c}(W) \rightarrow D^{b}_{c}(X)$ is an equivalence of derived categories. The essence of  
Theorem 2.3 in [BKR99] is the using 
of Serre duality and adjoint pairs. 

We remark that their argument  also shows that 
\[ \Psi : \Ext^{i}_{W} ( \mathcal{O}_{w_{1}}, \mathcal{O}_{w_{2}}) \rightarrow 
\Ext^{i}_{X} ( \mathcal{P}_{w_{1}}, \mathcal{P}_{w_{2}} ) \]
are isomorphisms for all $i$ (see [BKR99] p.18), from which one can  prove that $W$ is actually  a connected component of $M(X/Y)$. \newline
\newline \textbf{Step 4} \newline

The functor $\Psi: D^{b}_{c}(W) \rightarrow D^{b}_{c}(X)$ has a right adjoint $\Phi : D^{b}_{c}(X)  \rightarrow D^{b}_{c}(W)$, which is also a Fourier-Mukai type transform. The reader can see \cite{bkr:gnus} for an explicit formula. We show that $\Psi$ is fully faithful in this step. It suffices to show that $\Phi \circ \Psi \cong id$.

The composition functor 
$\Phi \circ \Psi$ is $\R         \pi _{2_{*}} ( \mathcal{Q}_{1} \bL \pi _{1}^{*}(-))$ where $\pi_{1}, \; \pi_{2}$ are the projections $W \times W \rightarrow W$ and $\mathcal{Q}_{1}$ is some object of $D(W \times W)$. It suffices to show that $\mathcal{Q}_{1}$ is quasi-isomorphic to $\mO_{\Delta_{W}}$.
We have $\cL i^{*}_{w} (\mathcal{Q}_{1})= \Phi \circ \Psi ( \mathcal{O}_{w})$.  
By an argument similar to the one given in Step 3, we have
$\Phi \circ \Psi ( \mathcal{O}_{w}) = \mathcal{O}_{w}$ for all $w$.
This shows that $\mathcal{Q}_{1}$ is actually the push-forward of a line bundle on $W$ to the diagonal $W \times W$. So $\Phi \circ \Psi$ is just twisting by $L$. To prove $\mathcal{Q}_{1}$ is quasi-isomorphic to 
 $\mO_{\Delta_{W}}$, 
it remains  to show $L$ is trivial.

There is a natural transform  $\varepsilon : id \rightarrow   \Phi \circ \Psi$, which gives a commutative diagram for every $w$: \newline
\xymatrix {
&&&&\mathcal{O}_{W} \ar[r]^{\varepsilon ( \mathcal{O}_{W})} \ar[d]_a & L \ar[d]^{L \otimes a} \\ 
&&&&\mathcal{O}_{w} \ar[r]^{ \varepsilon ( \mathcal{O}_{w})} & \mathcal{O}_{w} } \newline
where $a$ is non-zero. Since $\varepsilon$ is an isomorphism on the sub-category $D_{c}(W)$, it implies    $\epsilon ( \mathcal{O}_{W})$ is 
an isomorphism. This shows that $\mathcal{Q}_{1}$ is quasi-isomorphic 
to $\mO_{\Delta_{W}}$. \newline
\newline $\textbf{Step 5}$ \newline

By Lemma 2.1 in [BKR99], the statement that  the Fourier-Mukai type transform  $\Psi$ is an equivalence of derived categories follows from  the following statement
\[ \Phi (E) \cong 0 \Longrightarrow E \cong 0 \; \; \forall E \in D(X). \]
A proof of this statement can be found in Step 9 in \cite{bkr:gnus}.

  Harvard University, Department of Mathematics,  
                 One Oxford Street 
                 Cambridge, MA 02138 

 {\it E-mail address:} \tt{jcchen@math.harvard.edu}

\end{document}